\documentclass{amsart}
\usepackage{amssymb}
\usepackage{amscd}
\usepackage{latexsym}
\usepackage{amsfonts}
\usepackage{amsmath}
\newtheorem{proposition}{Proposition}
\newtheorem{theorem}{Theorem}
\newtheorem{definition}{Definition}
\newtheorem{conjecture}{Conjecture}
\newcommand{\eproof}{\begin{flushright} $\square$ \end{flushright}}

\newcommand{\Hom}{\mathop{\fam0 Hom}\nolimits}
\newcommand{\End}{\mathop{\fam0 End}\nolimits}

\newcommand{\ad}{\mathop{\fam0 ad}\nolimits}

\newcommand{\id}{\mathop{\fam0 Id}\nolimits}

\newcommand{\Aut}{\mathop{\fam0 Aut}\nolimits}

\newcommand{\Pic}{\mathop{\fam0 Pic}\nolimits}

\newcommand{\bC}{{\mathbb C}}

\newcommand{\C}{C}

\newcommand{\Z}{{\mathbb Z}}
\newcommand{\bZ}{\Z{}}

\newcommand{\ra}{\mathop{\fam0 \rightarrow}\nolimits}

\newcommand{\M}{{\mathcal M}}
\newcommand{\N}{{\mathbb N}}

\newcommand{\T}{ {\mathcal T}}
\newcommand{\fg}{ {\mathfrak g}}
\newcommand{\V}{ {\mathcal H}}
\renewcommand{\L}{{\mathcal L}}
\newcommand{\Sk}{{\mathcal S}}
\newcommand{\E}{ {\mathcal E}}

\renewcommand{\P}{ {\mathbb P}}

\newcommand{\s}{\sigma}
\newcommand{\Nabla}{{\mathbf {\hat \nabla}}}
\newcommand{\Nablae}{{\mathbf {\hat \nabla}}^e}

\newcommand{\Teim}{Teichm{\"u}ller }

\begin{document}
\title[The Nielsen-Thurstion classification via TQFT]{The Nielsen-Thurston classification of mapping classes is determined by
TQFT}
\author{J\o rgen Ellegaard Andersen}
\address{Department of Mathematics\\
        University of Aarhus\\
        DK-8000, Denmark}
\email{andersen@imf.au.dk}

\begin{abstract}
For each fixed $n\geq 2$ we show how the Nielsen-Thurston
classification of mapping classes for a closed surface of genus
$g\geq 2$ is determined by the sequence of quantum
$SU(n)$-representations $(\rho_k)_{k\in \N}$. That this is the
case is a consequence of the asymptotic faithfulness property
proved in \cite{A3}. We here provide explicit conditions on $(\rho_k(\phi))_{k\in
\N}$, which determines the Nielsen-Thurston type of any mapping
class $\phi$.
\end{abstract}

\maketitle

\section{Introduction}
The Nielsen-Thurston classification of mapping classes of compact
oriented surfaces splits mapping classes into three disjoint
types \cite{Th} (see also \cite{FLP} and \cite{BC}). We shall only be interested in the
closed surface case here, so we state the Nielsen-Thurston theorem
in this case.

Suppose $\Sigma$ is a closed oriented surface of genus $g\geq 2$ and let $\Gamma$ be the
mapping class group of $\Sigma$.
\begin{theorem}[Nielsen-Thurston]\label{NTclas}
A mapping class $\phi\in \Gamma$ has exactly one of the following
three properties.
\begin{itemize}
\item[1.] The mapping class $\phi$ is {\em finite order}, that is $\phi$
is a finite order element in $\Gamma$. This is equivalent to
$\phi$ having an automorphism of a Riemann surface as representative.
\item[2.] The mapping class $\phi$ is not finite order, but it is {\em reducible},
meaning there exists a simple closed curve on the surface, whose
non-trivial homotopy class is preserved by some power of $\phi$.
\item[3.] The mapping class $\phi$ is {\em Pseudo-Anosov}, meaning that
there exists $\lambda>1$, two transverse measured foliations $F^s$
and $F^u$ on $\Sigma$ and a diffeomorphism $f$ of $\Sigma$, which
represents $\phi$, such that
\[f_*(F^s) = \lambda^{-1} F^s \mbox{ and } f_*(F^u) = \lambda F^u.\]
\end{itemize}
\end{theorem}
In the Pseudo-Anosov case, $\lambda$ is uniquely determined by
$\phi$ and it is called the streching factor for $\phi$.

In the reducible case one continues the analysis of $\phi$, by
cutting $\Sigma$ along the preserved simple closed curve, to get a
mapping class of a surface with boundary. This mapping class is
then classified in terms of the Nielsen-Thurston classification of
mapping classes of surfaces with boundary. The upshot of this is
that there is a diffeomorphism $f$ of $\Sigma$, which represents
$\phi$ and preserve a system of simple closed curves, and when one
cuts the surface along these curves, $f$ induces a diffeomeorphism
of the resulting components, which on each piece is either finite
order or Pseudo-Anosov. See \cite{FLP} for further details
regarding this.

Fix an $n\geq 2$ and $d\in \Z/n\Z$. Let $\rho_k$ be the quantum
$SU(n)$ projective representations of the mapping class group $\Gamma$ at
level $k$ (with twist $d$ as described in \cite{A3}), which arises
from the Reshetikhin-Turaev $SU(n)$ Topological Quantum Field
Theory. These where first discussed by Witten in \cite{W1}. Then
they where rigorously constructed by Reshetikhin and Turaev in
\cite{RT1} and \cite{RT2}. Subsequently they where constructed
using skein theory in \cite{BHMV1}, \cite{BHMV2} and \cite{B1}.

The objective of this paper is to provide explicit conditions on the
endomorphisms $\rho_k(\phi)$ which determines which Nielsen-Thurston
type $\phi\in\Gamma$ belongs to.

Recall that we proved the asymptotic faithfulness property of the
sequence $\rho_k$ in \cite{A3}:
\begin{theorem}\label{MainA3} Assume that $n$ and $d$ are coprime or
  that $(n,d)=(2,0)$ when $g=2$.
Then we have that
\begin{equation*}
\bigcap_{k=1}^\infty \ker(\rho_k) =
\begin{cases}
\{1, H\} & g=2 \mbox{, }n=2 \mbox{ and } d=0 \\ \{1\}&
\mbox{otherwise}.
\end{cases}
\end{equation*}
where $H$ is the hyperelliptic involution.
\end{theorem}
The non coprime cases is covered in \cite{A4}\footnote{In
\cite{FWW} our argument from \cite{A3} was translated to the
BHMV-skein model. See also \cite{M2}.}.

Logically it follows from this Theorem, that the Nielsen-Thurston
classification of a mapping class $\phi$ is determined by
$(\rho_k(\phi))_{k\in\N}$, since this sequence determines $\phi$
itself. - However, we would here like to provide explicit conditions on
$(\rho_k(\phi))_{k\in\N}$, which separates the three
Nielsen-Thurston types.

The asymptotic faithfulness property gives us immediately the
following theorem.
\begin{theorem}\label{NTclasvTQFT}
For any mapping class $\phi\in \Gamma$ we have that there exist an
integer $M$ such that
\[\left( \rho_k(\phi)\right) ^M \in \bC \id\]
for all $k$ if and only if $\phi^M = 1$ (or $\phi^{2M} = 1$, in
case
 $(n,d) = (0,2)).$
\end{theorem}
This separates the finite order ones from the rest. In order to
separate the reducibles from the Pseudo-Anosov, we
recall the construction of the quantum representations and the
Toeplitz operator construction from \cite{A2} and \cite{A3}.

Let $p$ a point on $\Sigma$. Let $M$ be the moduli space of flat
$SU(n)$-connections on $\Sigma - p$ with holonomy $d\in \Z/n\Z
\cong Z_{SU(n)}$ around $p$. Assume that $n$ and $d$ are coprime
or that $(n,d)=(2,0)$ when $g=2$. Then $M$ is a smooth manifold.

By applying geometric quantization at level $k$ to the moduli
space $M$ one gets a vector bundle $\V^{(k)}$ over \Teim space
$\T$. The fiber of this bundle over a point $\sigma\in \T$ is
$\V^{(k)}_\sigma = H^0(M_{\sigma},\L_\s^k)$, where $M_{\s}$ is
$M$ equipped with a complex structure induced from $\sigma$ and
$\L_\s$ is an ample generator of the Picard group of $M_{\s}$.

The main result pertaining to this bundle $\V^{(k)}$ is that its
projectivization $\P(\V^{(k)})$ supports a natural flat
connection. This is a result proved independently by Axelrod,
Della Pietra and Witten \cite{ADW} and by Hitchin \cite{H}. Now,
since there is an action of the mapping class group $\Gamma$ of
$\Sigma$ on $\V^{(k)}$ covering its action on $\T$, which
preserves the flat connection in $\P(\V^{(k)})$, we get for each
$k$ a finite dimensional projective representation, say $\rho_k$,
of $\Gamma$, namely on the covariant constant sections of
$\P(\V^{(k)})$ over $\T$. This sequence of projective
representations $\rho_k$, $k\in {\mathbb N}$ is the {\em quantum
$SU(n)$ representations} of the mapping class group $\Gamma$.

For each $f\in C^\infty(M)$ and each point $\sigma\in \T$, we have
the {\em Toeplitz operator}
\[T^{(k)}_{f,\sigma} : H^0(M_{\sigma},\L_\s^k) \ra H^0(M_{\sigma},\L_\s^k)\]
which is given by
\[T^{(k)}_{f,\sigma} = \pi^{(k)}_\sigma(fs)\]
for all $s\in H^0(M_{\sigma},\L_\s^k)$, where $\pi^{(k)}_\sigma$
is the orthogonal projection onto $H^0(M_{\sigma},\L_\s^k)$.
We gets smooth section of
$\End(\V^{(k)})$ over $\T$
\[T_{f}^{(k)} \in C^\infty(\T,\End(\V^{(k)})) \]
by letting $T_{f}^{(k)}(\sigma) = T_{f,\sigma}^{(k)}$ (see
\cite{A3}).

The $L_2$-inner product on $\C^\infty(M,\L^k)$ induces an inner
product on $H^0(M_\sigma, \L_\sigma^k)$, which in turn induces the
operator norm $\|\cdot\|$ on $\End(H^0(M_\sigma, \L_\sigma^k))$.
Hence for any $A\in C^\infty(\T,\End(\V^{(k)}))$ we get the smooth
function $\|A\|$ on $\T$.

Suppose now $\gamma$ is a closed curve on the surface $\Sigma$.
Then we have the holonomy function $h_\gamma$ defined on $M$
associated to $\gamma$, given by taking the trace of the holonomy
around $\gamma$. Note that $h_\gamma$ only depends on the free
homotopy class of $\gamma$. Further $h_\gamma$ is constant if
$\gamma$ is nul-homotopic.

\begin{theorem}\label{NTClasvTQFT}
For any mapping class $\phi\in \Gamma$ and any
 homotopy class $\gamma$ of a simple closed curve on $\Sigma$
we have that $\phi$ is reducible along
$\gamma$,
i.e.
\[\phi(\gamma) = \gamma\]
if and only if $T^{(k)}_{h_\gamma}$ asymptotically commutes with $\rho_k(\phi)$:
$$ \lim_{k\ra \infty}
\|[\rho_k(\phi),T^{(k)}_{h_\gamma}]\| = 0.$$
\end{theorem}
Here the limit is pointwise convergence over $\T$. In fact, it
follows from the results in \cite{A3} that if the limit is zero at
some point in $\T$ then it holds for all points in $\T$ and in
fact the convergence is uniform on compact subsets of $\T$.

Using these two conditions, we see immediately how to determine
the Nielsen-Thurston classification of a mapping class:
\vskip.3cm
\begin{itemize}
\item[1.] Given a mapping class $\phi$, we first determine if it is
finite order or not by checking if $\rho_k(\phi)$ is finite order
bounded in $k$.
\item[2.] If not, we check through the non-trivial homotopy
classes $\gamma$
of simple closed curves to see if there exist
an integer $M$
with the property that $T^{(k)}_{h_\gamma}$ asymptotically commutes with $\rho_k(\phi^M)$:
$$ \lim_{k\ra \infty} ||[\rho_k(\phi^M),T^{(k)}_{h_\gamma}]|| = 0.$$ If so, $\phi$ is reducible.
\item[3.] If not, then $\phi$
is Pseudo-Anosov.
\end{itemize}
\vskip.3cm
We expect the following should be true.
\begin{conjecture}[AMU]
For a mapping class $\phi\in \Gamma$ we have that $\phi$
is Pseudo-Anosov or reducible with Pseudo-Anosov pieces if and
only if $\rho_k(\phi)$ is infinite order for large enough $k$.
\end{conjecture}
Evidence for this conjecture is provided by \cite{M} and
\cite{AMU}.

This paper is organized as follows.
In section \ref{sec2} we
recall the needed gauge theory setup. In the process we prove that
the free homotopy class of a simple closed curve on a surface is
determined by its holonomy function on the $SU(n)$ moduli space (for each $n\geq 2$).
In section \ref{sec4} we recall the main analytic estimate from
\cite{A3}, which states that the Toeplitz operators are
asymptotically flat with respect to Hitchin's connection.
The proof of Theorem \ref{NTClasvTQFT} is given in section
\ref{sec6}. In the last section we translate the formulation of the
results of this paper in the BHMV skein model for these TQFT's.

We would like to thank  Gregor Masbaum, Bill
Goldman and Peter Storm for discussions related to this work.

\section{The gauge theory construction of the quantum $SU(n)$ representations.}\label{sec2}
Let us now very briefly recall the construction of the quantum
$SU(n)$ representations. Only the details needed in this paper
will be given. We refer e.g. to \cite{H}, \cite{A3} and \cite{A5}
for further details. As in the introduction we let $\Sigma$ be a
closed oriented surface of genus $g\geq 2$ and $p\in \Sigma$. Let
$P$ be a principal $SU(n)$-bundle over $\Sigma$. Clearly, all such
$P$ are trivializable. As above let $d\in \Z/n\Z \cong Z_{SU(n)}$.
Throughout the rest of this paper we will assume that $n$ and $d$
are coprime, although in the case $g=2$ we also allow $(n,d) =
(2,0)$. Let $M$ be the moduli space of flat $SU(n)$-connections in
$P|_{\Sigma - p}$ with holonomy $d$ around $p$. We can identify
\[M = \Hom_d({\tilde \pi}_1(\Sigma), SU(n) )/ SU(n).\]
Here ${\tilde \pi}_1(\Sigma)$ is the universal central extension
\[0 \ra \bZ \ra {\tilde \pi}_1(\Sigma) \ra \pi_1(\Sigma) \ra 1\]
as discussed in \cite{H} and in \cite{AB} and $\Hom_d$ means the
space of homomorphisms from ${\tilde \pi}_1(\Sigma)$
to $SU(n)$ which send the image of $1\in \bZ$ in
${\tilde \pi}_1(\Sigma)$ to $d$ (see \cite{H}). When
$n$ and $d$ are coprime, $M$ is a compact smooth manifold of dimension $m = (n^2
-1)(g-1)$. In general, when $n$ and $d$ are not coprime $M$ is not
smooth, except in the case where $g=2$, $n=2$ and $d=0$. In
this case $M$ is in fact diffeomorphic to  $ {\mathbb C}\P^3$.
There is a natural homomorphism from the mapping class group
to the outer automorphisms of ${\tilde \pi}_1(\Sigma)$, hence
$\Gamma$ acts on $M$.

 We have that $M$ is a
component of the real slice in $\M$, the moduli space of flat
$SL(n,\C)$-connections on $\Sigma - p$, whose holonomy around $p$
has trace $n \exp(2 \pi \sqrt{-1} \frac{d}{n})$. This space is
connected and it contains the generalized \Teim space
$\tilde{\T}_p$ of $\Sigma - p$.

Suppose $\gamma$ is a closed curve on $\Sigma-p$. Then we have the
algebraic holonomy function $h_\gamma$ defined on the variety
$\M$. We are in particular interested in the restriction of
$h_\gamma$ to $M$, where it is a smooth function. We observe that
this restriction map is injective, since $M$ is a real slice in
$\M$, which is connected.
\begin{proposition}\label{Uniquec}
If $\gamma$ is a simple closed curve on $\Sigma-p$, then
the holonomy function $h_\gamma$ restricted to $M$
determines $\gamma$ up to free homotopy on $\Sigma
$.
\end{proposition}
\proof Suppose we have two simple closed curves $\gamma_1$ and
$\gamma_2$, such that $h_{\gamma_1} = h_{\gamma_2}$. Then if the
geometric intersection number of the homotopy classes of the two
curves is positive, we see that the two functions cannot possibly
be identical. This follows from considering the functions
restricted to the classical Teichm\"{u}ller space for the punctured
surface $\Sigma - p$. Consider a sequence of metrics where the
length of say $\gamma_1$ goes to zero. Then the length of
$\gamma_2$ most go to infinity. But the length of $\gamma_i$ is
determined by $h_{\gamma_i}$, so we get a contradiction.

Now suppose the geometric intersection number of the homotopy
classes of the two curves is zero. Then we can vary the holonomy
along the two curves completely independently, if they are not freely
homotopic. Hence they must be freely homotopic on $\Sigma$.
\eproof

Continuing towards the gauge theory definition of the quantum $SU(n)$ representations,
we now choose an invariant bilinear form $\{\cdot,\cdot\}$ on $\fg =
\mbox{Lie}(SU(n))$, normalized such that
$-\frac{1}{6}\{\vartheta\wedge [\vartheta\wedge\vartheta]\}$ is a
generator of the image of the integer cohomology in the real
cohomology in degree $3$ of $SU(n)$, where $\vartheta$ is the
$\fg$-valued Maurer-Cartan $1$-form on $SU(n)$.

This bilinear form induces a symplectic
form on $M$.  In fact
\[T_{[A]}M \cong H^1(\Sigma,d_A),\]
where $A$ is any flat connection in $P$ representing a point in $M$ and $d_A$ is the induced
covariant derivative in the associated adjoint bundle. Using this
identification, the symplectic form on $M$ is:
\[\omega(\varphi_1,\varphi_2) = \int_{\Sigma}\{\varphi_1\wedge\varphi_2\},\]
where $\varphi_i$ are $d_A$-closed $1$-forms on $\Sigma$ with
values in $\ad P$. See e.g. \cite{H} for further details on this.
The natural action of $\Gamma$ on $M$ is symplectic.

Let $\L$ be the Hermitian line bundle over $M$ and $\nabla$ the
compatible connection in $\L$ constructed by Freed in \cite{Fr}.
This is the content of Corollary 5.22, Proposition 5.24 and
equation (5.26) in \cite{Fr} (see also the work of Ramadas, Singer and Weitsman
\cite{RSW}). By Proposition 5.27 in \cite{Fr} we have that the
curvature of $\nabla$ is $\frac{\sqrt{-1}}{2\pi}\omega$. We will also use
the notation $\nabla$ for the induced connection in $\L^k$, where
$k$ is any integer.

By an almost identical construction, we can lift the action of
$\Gamma$ on $M$ to act on $\L$ such that the Hermitian connection
is preserved (See e.g. \cite{A1}). In fact, since $H^2(M, \bZ)\cong
\bZ$ and $H^1(M,\bZ) = 0$, it is clear that the action of $\Gamma$
leaves the isomorphism class of $(\L,\nabla)$ invariant, thus
alone from this one can conclude that a central extension of
$\Gamma$ acts on $(\L,\nabla)$ covering the $\Gamma$ action on
$M$. This is actually all we need in this paper, since we are only
interested in the projectivized action.

Let now $\sigma\in \T$ be a complex structure on $\Sigma$. Let us
review how $\sigma$ induces a complex structure on $M$ which is
compatible with the symplectic structure on this moduli space. The
complex structure $\s$ induces a $*$-operator on $1$-forms on $\Sigma$
and via Hodge theory we get that
\[H^1(\Sigma,d_A) \cong \ker(d_A + * d_A *).\]
On this kernel, consisting of the harmonic $1$-forms with values in
$\ad P$, the $*$-operator acts and its square is $-1$, hence we get an
almost complex structure on $M$ by letting $I= I_\sigma = -*$. It is a
classical result by Narasimhan and Seshadri (see \cite{NS1}),
that this actually makes $M$ a smooth
K{\"a}hler manifold, which as such, we denote $M_\sigma$. By
using the $(0,1)$ part of $\nabla$ in $\L$, we get an induced
holomorphic structure in the bundle $\L$. The resulting
holomorphic line bundle will be denoted
$\L_\s$. See also \cite{H} for further details on this.

From a more algebraic geometric point of view, we consider
the moduli space of S-equivalence classes of semi-stable bundles
of rank $n$ and determinant isomorphic to the line bundle ${\mathcal O}(d [p])$.
By using Mumford's Geometric
Invariant Theory, Narasimhan and Seshadri (see \cite{NS2}) showed
that this moduli space is a smooth complex algebraic projective variety
which is isomorphic as a
K{\"a}hler manifold to $M_\sigma$. Referring to \cite{DN} we
recall that
\begin{theorem}[Drezet \& Narasimhan]\label{DNTh1}
The Picard group of $M_\sigma$ is generated by the holomorphic line bundle
$\L_\s$ over $M_\sigma$
constructed above:
\[\Pic(M_\sigma) = \langle \L_\s \rangle. \]
\end{theorem}
\begin{definition}\label{Verlinde}
The bundle $\V^{(k)}$ over \Teim space is by definition the bundle
whose fiber over $\sigma\in \T$ is $H^0(M_\sigma,\L_\s^k)$, where
$k$ is a positive integer.
\end{definition}
There is a projective flat connection $\Nabla$ in the bundle
$\V^{(k)}$ over \Teim space, which is $\Gamma$-invariant due to
Axelrod, Della Pietra and Witten's and Hitchin's. We refer to
\cite{H} and \cite{ADW} for further details. See also \cite{A3}
and \cite{A5}.
Faltings has established this theorem in the case where one
replaces $SU(n)$ with a general semi-simple Lie group (see
\cite{Fal}). We further remark about genus $2$, that \cite{ADW} covers this case, but
\cite{H} excludes this case, however, the work of Van Geemen and De
Jong \cite{vGdJ} extends Hitchin's approach to the genus $2$ case.
\begin{definition}
The quantum $SU(n)$ projective representation $\rho_k$ at level $k$ of the
mapping class group is given by $\Gamma$'s action on the covariant
constant sections of $(\P(\V^{(k)}),\Nabla)$. We use the notation
$\P(Z^{(k)}(\Sigma))$ for this projective space of covariant
constant sections.
\end{definition}
We denote the induced flat connection in the endomorphism bundle
$\End(\V^{(k)})$ by $\Nablae$.

\begin{definition}
Let $\E^{(k)}(\Sigma)$ be the finite dimensional algebra of covariant
constant
sections of $(\End(\V^{(k)}), \Nablae)$ over $\T$ and let
\[\rho^e_{k} : \Gamma \ra \Aut(\E^{(k)})\]
be the corresponding representation.
\end{definition}

\section{Toeplitz operators and their asymptotic flatness}\label{sec4}

On $\C^\infty(M,\L^k)$ we have the $L_2$-inner product:
\[\langle s_1, s_2 \rangle = \frac{1}{m!}\int_M (s_1,s_2) \omega^m\]
where $s_1, s_2 \in \C^\infty(M,\L^k)$ and $(\cdot,\cdot)$ is the
fiberwise Hermitian structure in $\L^k$.

Now let $\sigma\in \T$. Then this $L_2$-inner product gives the
 orthogonal projection $$\pi^{(k)}_\sigma : \C^\infty(M,\L^k) \ra
 H^0(M_\sigma,\L_\sigma^k).$$
For each $f\in \C^\infty(M)$ consider the associated {\em Toeplitz
operator} $T_{f,\sigma}^{(k)}$ given as the composition of the
multiplication operator $$M_f : H^0(M_\sigma,\L_\sigma^k) \ra
\C^\infty(M,\L^k)$$ with the orthogonal projection:
\[T_{f,\sigma}^{(k)}(s) = \pi^{(k)}_\sigma(f s).\]
Then $T_{f,\sigma}^{(k)} \in \End(H^0(M_\sigma,\L_\sigma^k))$, and
we get a smooth section
\[T_{f}^{(k)} \in C^\infty(\T,\End(\V^{(k)})) \]
by letting $T_{f}^{(k)}(\sigma) = T_{f,\sigma}^{(k)}$ (see
\cite{A3}).

The $L_2$-inner product on $\C^\infty(M,\L^k)$ induces an inner
product on $H^0(M_\sigma, \L_\sigma^k)$, which in turn induces the
operator norm $\|\cdot\|$ on $\End(H^0(M_\sigma, \L_\sigma^k))$.

We need the following Theorem on Toeplitz operators
due to Bordemann, Meinrenken and Schlichenmaier (see \cite{BMS}, \cite{Sch},
\cite{Sch1} and \cite{Sch2}).
\begin{theorem}[Bordemann, Meinrenken and Schlichenmaier]\label{BMS1}
For any $f\in \C^\infty(M)$ we have that
\[\lim_{k\ra \infty}\|T_{f,\sigma}^{(k)}\| = \sup_{x\in M}|f(x)|.\]
\end{theorem}
Since the association of the sequence of Toeplitz operators
$T^k_{f,\sigma}$, $k\in \Z_+$ is linear in $f$, we see from this
Theorem, that this association is faithful.

Theorem 10 in \cite{A3} tells us that the Toeplitz operators are
asymptotically flat with respect to $\Nablae$:
\begin{theorem}\label{Asympflat}
Let $\sigma_0$ and $\sigma_1$ be two points in \Teim space and
$P_{\sigma_0,\sigma_1}$ be the parallel transport in the flat
bundle $(\End(\V^{(k)}),\Nablae)$ from $\sigma_0$ to $\sigma_1$. Then
\[\|P_{\sigma_0,\sigma_1}T_{f,\sigma_0}^{(k)} - T_{f,\sigma_1}^{(k)}\| = O(k^{-1}),\]
where $\|\cdot\|$ is the operator norm on $H^0(M_{\s_1},\L_{\s_1}^k)$.
\end{theorem}

\section{Proof of theorem \ref{NTClasvTQFT}}\label{sec6}

First we establish the following Proposition.
\begin{proposition}\label{Main1}
For any $\phi\in \Gamma$, $f\in C^\infty(M)$ and $\sigma\in \T$,
we have that $$ \lim_{k\ra \infty} \|T_{f,\sigma}^{(k)} - \rho_k(\phi)
T^{(k)}_{f,\sigma}\rho_k(\phi^{-1})\| = \lim_{k\ra \infty}\|T^{(k)}_{(f- f\circ\phi),
\sigma}\|.$$
\end{proposition}
\proof
Suppose we have a $\phi \in \Gamma$. Then $\phi$ induces a symplectomorphism of $M$
which we also just denote $\phi$ and we get the following commutative
diagram for any $f\in \C^\infty(M)$
\begin{equation*}
\begin{CD}
H^0(M_\sigma,\L_{\s}^k) @>\phi^*>> H^0(M_{\phi(\sigma)},\L_{\phi(\sigma)}^k)
@> P_{\phi(\sigma),\sigma}>> H^0(M_\sigma,\L_{\s}^k)\\
@V T^{(k)}_{f,\sigma} VV @V T^{(k)}_{f \circ \phi,\phi(\sigma)} VV
@VV{P_{\phi(\sigma),\sigma}T^{(k)}_{f \circ \phi,\phi(\sigma)}}V\\
H^0(M_\sigma,\L_{\s}^k) @>\phi^*>> H^0(M_{\phi(\sigma)},\L_{\phi(\sigma)}^k)
@> P_{\phi(\sigma),\sigma}>> H^0(M_\sigma,\L_{\s}^k),
\end{CD}
\end{equation*}
where $P_{\phi(\sigma),\sigma} :
H^0(M_{\phi(\sigma)},\L_{\phi(\sigma)}^k) \ra
H^0(M_{\sigma},\L_{\s}^k)$ on the horizontal arrows refer to
parallel transport in the bundle $\V^{(k)}$, whereas
$P_{\phi(\sigma),\sigma}$ refers to the parallel transport in the
endomorphism bundle $\End(\V^{(k)})$ in the last vertical arrow. By
the definition of $\rho_k$, we see that $$\rho_k(\phi)
T^{(k)}_{f,\sigma}\rho_k(\phi^{-1}) = P_{\phi(\sigma),\sigma}
T^{(k)}_{f\circ \phi,\phi(\sigma)}.$$ By Theorem \ref{Asympflat}
we get that
\[
\begin{split}
\lim_{k\ra \infty}\|T_{(f - f \circ \phi),\sigma}^{(k)}\| &=
\lim_{k\ra \infty}\|T_{f,\sigma}^{(k)} - T_{f\circ
\phi,\sigma}^{(k)}\|\\ & = \lim_{k\ra \infty}\| T_{f,\sigma}^{(k)}
- P_{\phi(\sigma),\sigma} T^{(k)}_{f\circ \phi,\phi(\sigma)}\| \\
& = \lim_{k\ra \infty}\| T_{f,\sigma}^{(k)} - \rho_k(\phi)
T^{(k)}_{f,\sigma}\rho_k(\phi^{-1})\|.
\end{split}\]
\eproof
\proof[Proof of Theorem \ref{NTClasvTQFT}] Let $\phi\in \Gamma$
and assume that there exists $\gamma$ non nul-homotopic simple
closed curve on $\Sigma$ and $\sigma\in \T$, such that
$$
\lim_{k\ra \infty} \|[\rho_k(\phi),T^{(k)}_{h_\gamma, \sigma}]\| =
0.$$
But
\[\begin{split}\| T_{h_\gamma,\sigma}^{(k)} - \rho_k(\phi^{-1})
T^{(k)}_{h_\gamma,\sigma}\rho_k(\phi)\| & =  \|\rho_k(\phi^{-1}) \rho_k(\phi)(
T_{h_\gamma,\sigma}^{(k)} - \rho_k(\phi^{-1})
T^{(k)}_{h_\gamma,\sigma}\rho_k(\phi))\|\\
& \leq \|\rho_k(\phi^{-1})\|\|[\rho_k(\phi),T^{(k)}_{h_\gamma, \sigma}]\|.
\end{split}\]
Lemma 1 and Proposition 2 in \cite{A3} gives
a uniform bound on $\|\rho_k(\phi^{-1})\|$.
Hence Proposition \ref{Main1} implies that
\[\lim_{k\ra \infty} \|T^{(k)}_{(h_\gamma - h_\gamma \circ \phi^{-1}), \sigma}\| =
0.\] Then by Bordemann, Meinrenken and Schlichenmaier's Theorem
\ref{BMS1}, we must have that $h_\gamma = h_\gamma\circ \phi^{-1}$. But
so by Proposition \ref{Uniquec}, we must have that $\phi$
preserved the homotopy class of $\gamma$.

The converse follows directly from Proposition \ref{Main1}, since $h_\gamma = h_\gamma\circ\phi$,
if $\phi(\gamma) = \gamma$.
\eproof

\section{The formulation in the BHMV-skein model}

In this section we provide a translation of our conditions for separating
the Nielsen-Thurston types of mapping classes to the
BHMV-skein model. The determination of the finite order elements, Theorem \ref{NTclasvTQFT} needs
no translation, since it follows directly from the asymptotic faithfulness
property. We just need to translate the condition in Theorem
\ref{NTClasvTQFT} which determines when a mapping class is
reducible.

Let $V_p$ be the TQFT defined by Blanchet,
Habegger, Masbaum and Vogel in \cite{BHMV1} and \cite{BHMV2},
where $p = 2 r$, $r$ being an integer and we choose the $2p$ root of $1$
to be $A = e^{i \pi/2r}$. In particular there is also a
projective representation
\[V_p : \Gamma \ra \Aut(\P(V_p(\Sigma)))\]
for each $p$.

Let $\Sk(\Sigma)$ be the free $\bC$-vector space generated by
isotopy classes of one-dimensional sub-manifolds of $\Sigma$. We note
that $\Sk(\Sigma)$ is a subspace of the free $\bC$-vector space generated by
links in $\Sigma\times [0,1]$.

Since $V_p(\partial (\Sigma\times [0,1])) =
\End(V_p(\Sigma))$, we get by the very construction of the
BHMV-skein
model of this TQFT a sequence of linear maps
\[V_p : \Sk(\Sigma) \ra \End(V_p(\Sigma)).\]

\begin{theorem}\label{ReducibleBHMV}
For any mapping class $\phi\in \Gamma$ and any non-trivial
 homotopy class $\gamma$ of a simple closed curve on $\Sigma$
we have that $\phi$ is reducible along
$\gamma$,
i.e.
\[\phi(\gamma) = \gamma\]
if and only if
$$[V_p(\phi),V_p(\gamma)] = 0$$
for all $p=2r$.
\end{theorem}

Our joint work with K. Ueno \cite{AU1}, \cite{AU2} and \cite{AU3},
combined with the work of Laszlo \cite{La1}
gives the following result for the $(n,d) = (2,0)$ theory:
\begin{theorem}\label{projiso}
There is a projective linear isomorphism of representations of $\Gamma$
\[I_k :  \P(V_{2k+4}(\Sigma)) \ra \P(Z^{(k)}(\Sigma)).\]
\end{theorem}

The projective linear isomorphism $I_k$ induces a  linear
isomorphism of representations of $\Gamma$
\[I^e_k :  \End(V_{2k+4}(\Sigma)) \ra  \E^{(k)}(\Sigma).\]

In \cite{A6} we prove that
\begin{theorem} \label{Toeplitzgamma}
For any non-trivial homotopy class $\gamma$ of a simple closed curve on $\Sigma$
we have that
\[\lim_{k\ra \infty} \|I_k^e(V_{2k+4}(\gamma)) - T^{(k)}_{h_\gamma}\| = 0.\]
\end{theorem}

\proof[Proof of Theorem \ref{ReducibleBHMV}]
Let $\phi\in\Gamma$ be a mapping class. Assume
there exists
$\gamma$ non nul-homotopic simple closed curve on $\Sigma$  such that
$$[V_{2k+4}(\phi),V_{2k+4}(\gamma)] = 0$$
for all $k\in\N$. Then it follows immediately that
\[\lim_{k\ra \infty} \|[I_k^e(V_{2k+4}(\phi)),T^{(k)}_{h_\gamma}]\| = 0\]
by Theorem \ref{Toeplitzgamma}.
But since $I_k^e(V_{2k+4}(\phi)) = \rho_k(\phi)$ by Theorem \ref{projiso}, we get by
Theorem \ref{NTClasvTQFT} that $\phi$ is reducible.
The converse statement is trivial.
\eproof

After the completion of this work we where made aware of
\cite{MN}. One can translate our proof of Theorem
\ref{NTClasvTQFT} into a BHMV-skein model proof of Theorem \ref{ReducibleBHMV}
using the result of \cite{MN}.

\end{document}